\documentclass{amsart}
\usepackage{amssymb}
\newtheorem{thm}{Theorem}[section]

\newtheorem{lem}[thm]{Lemma}

\theoremstyle{definition}

\theoremstyle{remark}

\numberwithin{equation}{section}

\newcommand{\eps}{\varepsilon}

\newcommand{\rt}{\rightarrow}

\newcommand{\C}{{\mathbb C}}
\newcommand{\B}{{\mathbb B}^n}
\newcommand{\pa}{\partial M}
\newcommand{\om}{\overline M}
\newcommand{\te}{\widetilde {exp}}

\begin{document}

\title[nonpositively curved K\"ahler manifolds]{A class of nonpositively curved
K\"ahler manifolds \\ biholomorphic to the unit ball in $\C^n$}
\author{harish seshadri}
\address{department of mathematics,
Indian Institute of Science, Bangalore 560012, India}
\email{harish@math.iisc.ernet.in}
\author{Kaushal Verma}
\address{department of mathematics,
Indian Institute of Science, Bangalore 560012, India}
\email{harish@math.iisc.ernet.in}


\begin{abstract}
Let $(M,g)$ be a simply connected complete K\"ahler manifold with
nonpositive sectional curvature. Assume that $g$ has constant
negative holomorphic sectional curvature outside a compact set.
We prove that $M$ is then biholomorphic to the unit ball in $\C
^n$, where $dim_{\C} M=n$.
\end{abstract}
\maketitle

 \noindent {\bf R\'esum\'e}. Soit $(M,g)$ une vari\'et\'e k\"ahl\'erienne compl\`ete
et simplement connexe \`a courbure sectionnelle non-positive.
Supposons que $g$ ait courbure sectionnelle holomorphe constante
et n\'egative en delors d'un compact. On d\'emontre que $M$ est
biholomorphe \`a une boule dans $\C^n$, o\`u \ $dim_{\C}M=n$.

\section{Introduction}
An important issue in complex differential geometry is to
understand the relationship between the curvature of a K\"ahler
manifold and the underlying complex structure. The simplest
theorem along these lines is the classical theorem of Cartan that
a simply connected complete K\"ahler of constant holomorphic
sectional curvature is holomorphically isometric to $ {\mathbb
C}P^n$, $\B$ or $\C ^n$ (where $\B$ is the open unit ball in
$\C^n$) depending on whether the curvature is positive, negative
or zero. Here the metrics on $ {\mathbb C}P^n$, $\B$ or $\C ^n$
are the Fubini-Study metric, the Bergman metric and the flat
metric, respectively.

A far deeper theorem is that of Siu and Yau ~\cite{sy} which
states that a complete simply connected nonpositively curved
K\"ahler manifold of faster than quadratic curvature decay has to
be biholomorphic to $\C^n$. An analogue of this theorem for
characterizing the ball in $\C^n$ is not known, to the best of
our knowledge. As a first step in this direction we prove the
following theorem, which can also be regarded as perturbed
version of Cartan's theorem stated above, at least in the
negative case.

\begin{thm}\label{mai}
Let $(M,g)$ be a simply connected complete K\"ahler manifold with
nonpositive sectional curvature. If $g$ has constant negative
holomorphic sectional curvature outside a compact set, $M$ is
biholomorphic to the unit ball in $\C ^n$, where $dim _{\C}M=n$.
\end{thm}
It is natural to ask if the theorem is true if we only have that
the holomorphic sectional curvatures converge to $-1$ as $r \rt
\infty$, where $r$ is the distance from a fixed point in $M$.
However, the following class of examples show that the theorem
then fails: If $g$ is the Bergman metric of a strongly
pseudoconvex domain $\Omega$ in $\C^n$, then the holomorphic
sectional curvatures of $g$ approach $-1$ as one approaches
$\partial \Omega$. Moreover, by the results of ~\cite{gk}, if
$\Omega$ is a small enough perturbation of $\B$, i.e., $\partial
\Omega$ is a $C^\infty$ small perturbation of $\partial \B$, then
$g$ has negative sectional curvature. By Chern-Moser
theory,``many" of these perturbations are not biholomorphic to
$\B$. Hence one should impose a specified rate of convergence of
holomorphic sectional curvatures to $-1$ in order to obtain a
theorem similar to that of Siu-Yau in the negative case. However,
it is not clear what this rate of convergence should be.

Finally, we note that the hypotheses of the Siu-Yau theorem are
strong enough to guarantee that the K\"ahler manifold is actually
holomorphically {\it isometric} to $\C^n$ with the flat metric.
In fact, R. Greene and H. Wu proved that a {\it Riemannian}
manifold with the same curvature hypotheses has to be flat
~\cite{gw}. In our case, however, we can always perturb the
Bergman metric of $\B$ on a compact set and satisfy our
hypotheses.

Roughly, the proof of Theorem \ref{mai} proceeds as follows:
Suppose that $M$ has constant holomorphic curvature outside a
compact set $K$. As in the proof of Cartan's theorem, one can use
the exponential map to construct holomorphic maps to $\B$ on
``pieces" of $M \setminus K$. The difficulty here is that even
though these maps can be chosen to patch up to a give single
holomorphic map from $M \setminus K$ to $\B$, this map may not be
injective. We avoid this difficulty by working with $\partial M$,
the asymptotic boundary of $M$. More precisely, we use the
holomorphic maps above to define a spherical CR-structure on
$\pa$. Since $\pa$ is simply connected, one gets a global
CR-diffeomorphism to $S^{2n-1}$. One then notes that since $M$ is
Stein, we can extend this diffeomorphism to $M$ by Hartogs'
theorem.

\section{Proof}
For the rest of this paper, $M$ will denote a simply-connected,
complete K\"ahler manifold with nonpositive sectional curvature
and constant holomorphic sectional curvature $-1$ outside a
compact set. $\pa$ will denote its asymptotic boundary. There is
a natural topology, described in the proof below, on \ $\overline
M:= M \cup \pa$ \ which makes it a compact topological
manifold-with-boundary.


The main theorem is proved by first proving the following
proposition. In what follows $S^{2n-1}$ is the unit sphere in
$\C^n$ with the induced CR-structure.


\begin{lem}\label{cr}
 $\overline M$ can be given the structure of a smooth compact manifold-with-boundary such that $\pa$ admits a ``natural" CR-structure
 which makes it CR-diffeomorphic to $S^{2n-1}$.



\end{lem}

Before beginning the proof, we recall certain general
constructions on nonpositively curved manifolds: \vspace{2mm}

First, we define the {\it ``modified"} exponential map. Let \ $V$
\ be a complex vector space with a Hermitian inner product
 \ $h$ \ and let \ $B= \{ x \in V: \Vert x \Vert < 1 \}$ \ denote the open
 unit ball in \ $V$. Define the homeomorphism \ $\phi:V \rt B$ \
 by \ $\phi(x)= (1-e^{- \Vert x
 \Vert}) \ \frac{x}{\Vert x \Vert}$. Note that \ $\phi$ \ is a diffeomorphism on \ $V \setminus
 \{0\}$.
 When \ $V=T_pM$ \ and \ $h=g_p$, we will use the notation
 \ $\phi_p$. For any \ $p \in M$, define the {\it modified exponential
 map} \ \
 $\te _p : B_p \rt M$ \ by  \ $ \te_p= exp_p \circ \ \phi_p^{-1}.$
 \vspace{2mm}

Next, let \ $\partial M \ = \{$ equivalence classes of geodesic
rays in $M \}$, where geodesics \ $c_1, c_2:[0, \infty) \rt M$ \
are
  equivalent if there is a constant \ $a < \infty$ \ such that \ $d(c_1(t),c_2(t)) <a$ \ for all \ $t \ge
  0$. \ $\pa$ \ is usually referred to as the {\it asymptotic boundary} of
 \ $M$. We endow \ $\overline M = M \cup \pa$ \ with the ``cone" topology. This is the topology generated by open sets in \ $M$ \ and
 ``cones", corresponding to \ $x \in M$, $z \in \pa$ \ and \ $\eps >0$, defined by
 \ $$C_x(z, \eps):=\{ y \in \overline M \ \vert \ y \neq x \ \ {\rm and}\ \ <_x(z,y) < \eps \}.$$
Here the angle \ $<_x(z,y) \ := \quad <(c_1^{'}(0),c_2^{'}(0))$ \
where \ $c_1$ \ and \ $c_2$ \ are geodesics joining  $x$  with
$z$  and
 $y$ (see ~\cite{gro} for details).
 For any \ $p \in M$, \ $\te_p$ \ extends to a homeomorphism, which we continue
 to denote by the same symbol, from \ $\overline B_p$ \ to \ $\overline
 M$.
 \vspace{2mm}





 Now we come to the proof of
Lemma \ref{cr}.
\begin{proof}
Suppose \ $M$ \ has constant holomorphic curvature $-1$ outside a
compact set \ $K$. Fix \ $o \in M$. Choose \ $R$ large so that \
$d(o,x) <R$ \ for any \ $x \in K$. \vspace{1mm}

If \ $p \in \pa$, then there is an unit-speed geodesic \ $\gamma_p
 :[0, \infty) \rt \om$ \ with \ $\gamma _p(0)=o$ \ and \ $lim _{t \rt
 \infty} \gamma_p (t)=p$.
Let
$$x(p):=\gamma _p(R).$$
 We observe that $C_{x(p)}(p, \frac {\pi}{4})
\cap K = \emptyset$. This is because $d(o,x)>R$ for any $x \in
C_{x(p)}(p, \frac {\pi}{4})$ by Toponogov's Comparison Theorem for
geodesic triangles in nonpositively curved manifolds (see
~\cite{gro}, Page 5). Hence, by our choice of $R$, $g$ has
constant holomorphic sectional curvature in the interior of
$C_{x(p)}(p, \frac {\pi}{4})$. \vspace{1mm}

 Choose \ $p_1,..,p_k \in \pa$ \ so that if
 $$ U_i := C_{x(p_i)}({p_i}, \frac {\pi}{4}) ,$$
 then \ $ {U}_{1} \cap \pa,.., {U}_{k} \cap
 \pa$ \
cover \ $\pa$.

For \ $i=1,..,k$, choose linear isometries \ $L_i: T_{\gamma_{p_i}
(R)} \rt T_0 \B$. We then get maps
$$f_{i}:=  \te_0 \circ L_i \circ \te_{x(p_i)}^{-1} :{U}_{i} \rt  \overline
\B.$$ These maps are homeomorphisms onto their images and we
declare these to be charts on \ $\pa \subset \overline M$. In
order to check that the transition functions are smooth, let us
observe the following:

First, it is easily checked that \ $ f_{i} \vert_{U_i}=  exp_0 \
\circ L_i \ \circ exp_{p_i}^{-1}$. Recall that our metric is
locally symmetric in the interior of \ $U_{i}$. Hence by the
Cartan-Ambrose-Hicks Theorem (cf. ~\cite{ce}), $f_{i}$ \ is a
holomorphic local isometry there.

Next, Toponogov's Comparison Theorem implies that for any $p \in
\pa$, $U_p$ is geodesically convex. Also, it is clear from the
definition that if $q \in U_p$, then the geodesic ray starting at
$x(p)$ and passing through $q$ lies in $U_q$. Combining these
observations, we see that if $U_i \cap U_j \neq \emptyset$, then
\begin{equation}\label{con}
U_i \cap U_j \ \ is \ \ connected \ \ and \ \ U_i \cap U_j \cap
\pa \neq \emptyset.
\end{equation}

Now the transition function $f_i \circ f_j^{-1}$ is a holomorphic
isometry (for the restriction of the Bergman metric of $\B$) from
$f_j(U_i \cap U_j) \cap \B$ to $f_i(U_i \cap U_j) \cap \B$. Since
$f_j(U_i \cap U_j)$ is connected by (\ref{con}), such a mapping
has to be the restriction of a global automorphism of $\B$. In
particular, the mapping is smooth up to the boundary, i.e. $f_i
\circ f_j^{-1}: f_j(U_i \cap U_j) \rt f_i(U_i \cap U_j)$ is
smooth.
This gives us the required smooth structure on $\overline M$.

Also, it is clear that the charts \ $(U_{i} \cap \pa, f_{i})$ \
define a CR-structure on $\pa$, since the transition functions
will be local CR-diffeomorphisms of \ $S^{2n-1}$. Moreover, by
definition this CR-structure is locally spherical. Since  $\pa$
is compact and simply connected, the results of ~\cite{bur}
(basically a developing map argument) imply that there is a
global diffeomorphism \ $\psi$ \ from \ $\pa$ \ to \ $S^{2n-1}$.
This proves the lemma. \vspace{1mm}
\end{proof}

We continue with the proof of the main theorem, using the
notation in the proof of the lemma. Let us note that by composing
with holomorphic automorphisms of \ $\B$, if necessary, we can
assume that \
\begin{equation}\label{sam}
f_{i}\vert_{U_{i} \cap \pa} = \psi.
\end{equation}
By (\ref{con}), (\ref{sam}) and unique continuation, if \ $U_i
\cap U_j \neq \emptyset$ \ then \ $f_i \circ f_j^{-1} = \ id: \
f_j (U_i \cap U_j) \rt f_i(U_i \cap U_j)$, \ since \ $f_i \circ
f_j^{-1}= \ id$ \ on \ $f_j(U_i \cap U_j) \cap \partial \B$.

 Hence if \ $U_i \cap U_j \neq \emptyset$, \ then \ $f_i=f_j$ \ on \ $U_i \cap U_j$. Therefore, the $f_i$ patch up
to give a
smooth mapping \ $  F:U \rt \B$ \ on the open neighbourhood \ $U=
U_1 \cup ... \cup U_k $ \ of \ $\pa$ \ and $F$ is holomorphic on \
$U \cap M$. Since \ $F|_{\pa}=\psi$ \ is injective and since $F$
is a local diffeomorphism, we can choose a neighbourhood $V
\subset U$ of $\pa$ such that $F|_V$ is injective.

To extend $F$ to $M$, we recall Wu's theorem ~\cite{wu} that {\it
a simply connected complete K\"ahler manifold of nonpositive
sectional curvature is Stein}. Combining this with the fact that
 $M \setminus U$  is compact, we conclude that  $ F$
extends to all of  $M$  by Hartogs' theorem on Stein manifolds.
By the maximum principle, \ $F(M) \subset \B$.

To construct the inverse of $F$, let \ $G=F|_V^{-1}$. \ $G$ is
smooth map defined on the neighbourhood $F(V)$ of \ $\partial \B$,
which is holomorphic in \ $F(V) \cap \B$.
Since $M$ is Stein, $M$ is an embedded submanifold of some \
${\mathbb C}^N$. Again by Hartogs theorem and the maximum
principle, \ $G: F(V) \cap \B \rt V \subset \C^N$ \ extends to a
smooth map \ $G: \overline \B  \rt \overline M$ \ which is
holomorphic in $\B$.

Finally, by unique continuation, \ $F \circ G =id_ { \overline
\B}$ \ and \ $G \circ F=id _{\overline M}$ \hfill Q.E.D.
\vspace{2mm}

\noindent {\bf Remark}: A CR-structure on the boundary of a
nonpositively curved K\"ahler manifold is shown to exist under
hypotheses more general than ours in ~\cite{bla}. \vspace{2mm}

\noindent We end with the following \vspace{2mm}

\noindent {\bf Question}: Let $(M^n ,g)$ be a simply-connected
complete K\"ahler manifold of nonpositive curvature. If $g$ is
locally symmetric outside a compact set, is $M$ biholomorphic to
$\Omega \times \C^{n-k}$, \ where $\Omega$ is a bounded symmetric
domain in $\C^k$, for some $k$ ?


\end{document}